# Finite Element Analysis of Shear Lag Effect in Long-Span Single-Box Continuous Rigid Bridges

Zhaokun Shen [1], Chengtao Zhou [2], Dianhao Li [1], Yun Meng [3], Lvye Zhou [1], Ligui Yang [4], Chengmao He [4], Shaorui Wang [4] and Shuangshuang Jin [4]

1   Guizhou Qiantong Engineering Technology Co., LTD,Guiyang 550000, China
2   Guizhou Province transportation comprehensive administrative law enforcement directly under the detachmen, Guiyang 550000, China
3   Guizhou Qiantong Anda Engineering Consulting Co., LTD, Guiyang 550000, China
4   School of Civil Engineering, Chongqing Jiaotong University, Chongqing 400074, China
*   Correspondence: chengtaozhou@163.com

**Abstract:** As the span and width of continuous rigid bridges increase, the complexity of the spatial forces acting on these structures also grows, challenging traditional design methods. Primary beam theory often fails to accurately predict the stresses in the bridge girder, leading to potential overestimation of the ultimate capacity of these bridges. This study addresses this gap by developing a detailed finite element (FE) model of a continuous rigid bridge using ABAQUS, which accounts for the complex 3D geometry, construction procedures, and nonlinear material interactions. A loading test on a 210-meter continuous rigid bridge is performed to validate the model, with the measured strain and deflection data closely matching the FE simulation results. The study also examines stress distributions under constant and live loads, as well as the shear lag coefficients of the main girder. Through parametric analysis, we explore the effects of varying the width-to-span and height-to-width ratios. The results reveal that both a wider box girder and a larger span significantly amplify the shear lag effect in the bridge girder. The findings enhance the understanding of stress distribution in large-scale rigid bridges and provide critical insights for more accurate design and assessment of such structures.

**Keywords:** continuous rigid bridge; large span wide box girder; finite element analysis; shear lag effect

## 1. Background

Large-span continuous rigid bridges are critical for crossing challenging terrains, such as small canyons and mountains, and play a key role in transportation infrastructure. However, as the span and width of these bridges increase, the complexity of the forces acting on them intensifies. In particular, uneven shear deformation in the box girder wings results in the distribution of positive stresses along the cross-section in a curved manner, complicating the stress characteristics of the structure. Although these bridges offer excellent overall performance, the initial designs of large-span wide box girder bridges, particularly those exceeding 200 meters, often fail to account for the spatial shear lag effect that can lead to the non-uniform compressive stress at the web-to-flange junction. As shown in Figure .1, it not only brings about cracking in the girder over time, but also the long-term deflection of the girder, complicating maintenance and potentially disrupting the normal operation of the bridge. Therefore, understanding the force mechanisms of large-span, wide-box girders is crucial to ensuring the long-term performance and reliability of such bridges.

The shear lag effect, a mechanical phenomenon that negatively impacts the load-bearing capacity of box-type components, can lead to a significantly uneven stress distribution across the cross-section of flexural and shear-loaded member [1]. The shear flow transmitted from the web to the flange plate leads to this lag, which prevents the flange from maintaining the assumption of a flat cross-section, a key principle in traditional beam theory [2]. As the span and aspect ratio of box girders increase, the shear lag effect becomes more pronounced due to the intensifying spatial forces and the non-uniform distribution of stresses across the cross-section [3].

Theoretical studies have been conducted on the force transfer mechanisms of large-span continuous rigid main girders[4-7], including the effective distribution-width method[8], the simulated rod method[9], the three-rod simulation method[10], the additional deflection method[11], the energy variational method[12], the Timoshenko beam-energy variational combination method[13], the orthotropic anisotropic plate method[14-15], the plate and shell theory method[16], and finite segment methods[17]. Despite the studies on the shear flow and force transfer mechanisms in large-span continuous girders [18], the consideration of boundary conditions is generally simplified for convenience. Most traditional approaches, such as elastic beam theory and simplified shear distribution methods, fail to adequately address the complexities of spatial shear effects [19]. As large cantilever single box single chamber wide girder bridges continue to develop, the complexity of spatial forces increases. Along with the span of continuous rigid bridge girder, the cross-sectional height is reduced, which can result in more severe shear lag effect and increasing stress nonuniformity.

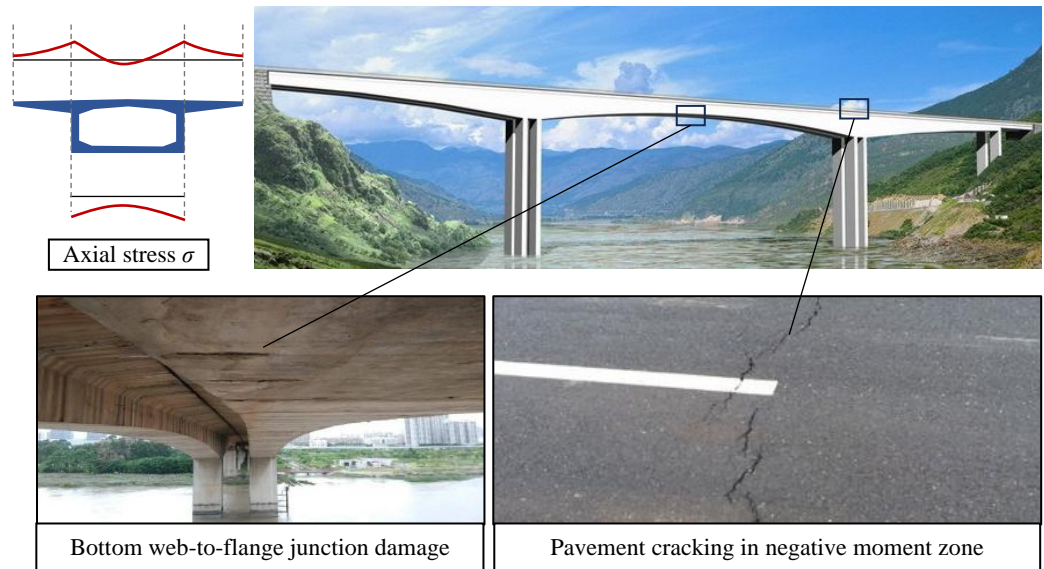

**Figure 1.** Diseases of concrete continuous rigid bridges

The failure to consider the stress nonuniformity caused by shear lag effect in the early design of continuous concrete rigid bridges has led to a number of practical issues. In many cases, these bridges exhibit excessive deflection and cracking, particularly as they age and undergo heavy loading. The issues have led to significant costs in bridge maintenance and repair, as well as disruptions to the normal operation of transportation networks. The increasing complexity of spatial inner forces in large-span bridges, along with the growing traffic volume and heavier loads, highlights the need for more accurate and comprehensive design methods. Recently, the finite element (FE) analysis has been employed in recent years to better understand the shear lag effect in box girders [20-22]. Yet, to better predict the performance and durability of these structures, the dynamic loading, nonlinear interactions, and material plasticity should be considered in future FE studies.

In order to accurately analyze the force state and spatial effects of large-span, wide-box and large-cantilever continuous rigid-frame bridges, this paper establishes a refined finite element model using ABAQUS to represent the complex spatial structure, and load testing on an actual bridge has been conducted to verify the accuracy of the simulation results. The study also investigates the stress distribution patterns under constant and live loads, as well as the distribution of the shear lag coefficient in the main girder. The results may help to improve the accuracy of existing design codes and contribute to the ongoing development of more advanced computational models for the analysis of complex bridge structures.

## 2. A typical long-span concrete continuous rigid frame bridge

In this paper, a verified finite element analysis model is developed based on a prestressed concrete continuous rigid frame bridge constructed in 2006 (Figure 2). The bridge consists of three segments with spans of 122 m, 210 m, and 122 m, featuring a single-cell box girder. The top plate has a width of 22.5 m, while the bottom plate is 11 m wide. The minimum cross-sectional height at the mid-span and in the two side spans is 3.5 m. The root height of the box girder is 12.5 m, and the top plate thickness ranges from 30 to 50 cm. The cantilever flange at the end has a thickness of 20 cm, and the web near the bridge pier is 120 cm thick, while the remaining portion of the web varies in thickness from 50 cm to 70 cm. The bottom plate thickness of the girder ranges from 32 cm to 150 cm. Each T-frame girder of the bridge is equipped with 174 prestressed steel strands, 42 of which are located at the mid-span. Among them, 40 strands pass through the bottom plate, while 2 strands pass through the top plate. The tensile control stress is set at 0.75 times the standard tensile strength (0.75fpk), equivalent to 1395 N/mm². The top plate and web of the box girder are reinforced with 22 strands of $\phi_s15.24$ steel. The spacing of the transverse prestressed steel strands within the top plate varies from 70 cm to 75 cm and consists of 4 strands of $\phi_s15.24$ steel. For vertical prestressed reinforcement, $\phi32$ high-strength rolled rebars are used, delivering a tensile force of 54.2 t.

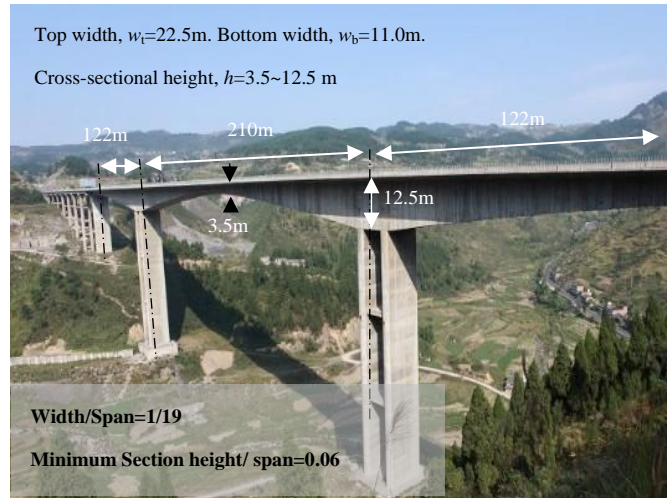

**Figure 2.** Geometry of bridge

## 3. Finite element model

The finite element of the bridge is established according to Figure 3, illustrating the reinforcement mesh and prestressing bundles embedded within the concrete box girder and abutment entities. The diaphragm of the box girder is connected to the inner wall of the box girder using the "Tie" connection. The edge of the main girder at the side span is articulated vertically, and a spring bearing simulates the rubber bearing laterally, with a single shear stiffness of $2\times10^4$ kN/m. In addition, the bottom of the bridge pier is rigidly connected.

### 3.1 Elements

The solid elements (C3D8R) were generated using the 'Sweep' method to represent concrete girders and piers. Truss elements (T3D2) were utilized for the prestressed steel strands, while rebar layer elements (SFM3D4) were employed for reinforcing bars. Additionally, shell elements (S4R) were used for the flange slabs and web thickening layers.

### 3.2 Materials

The compressive strength of concrete is 28.5 N/mm², while its tensile strength is 2.45 N/mm², accompanied by an elastic modulus of $3.5\times10^4$ N/mm². The behavior of concrete is simulated using the concrete damage plastic (CDP) model [23], which adopts the uniaxial stress-strain model outlined in Appendix C.2.3 of the Code for the Design of Concrete Structures GB50010-2010 [24]. Key parameters include the concrete expansion angle ($\psi$), tensile/compressive ratio (Kc), flow potential skewness, initial equivalent biaxial compressive yield stress, and the initial uniaxial compressive yield stress ratio, which are set to 36°, 2/3, 0.1, 1.1667, and 0.0005 [25], respectively.

The steel plates were designed by using Chinese Q345 steel, and the prestressing strands of the main girder comprised 22$\phi$s15.24 high-strength low-relaxation strands, which have a standard strength of 1860 MPa and a modulus of elasticity of $1.95\times10^5$ MPa. The tensioning control stress was set at 1395 MPa ($0.75\sigma_{con}$). The yield strength of the light round steel bar is 240 MPa, whereas the yield strength of Grade II steel bar is 340 MPa. The stress-strain relationship model for steel follows the ideal elastic-plastic bifurcation model outlined in Appendix C.1.2 of the Code for Structural Design of Concrete Structures GB50010-2010[24].

### 3.3 Loads

The permanent load refers to the self-weight of each component, primarily consisting of the concrete of the bridge box, prestressed tendons, transverse and longitudinal steel bars, diaphragms, and the bridge deck pavement. The self-weight of the concrete and steel layer units is applied through 'Gravity', with a gravitational acceleration set at 9.8 m/s². Due to the use of truss units in the prestressed beams, their own weight can easily result in non-convergence when applied via gravity; therefore, an equivalent load is applied to the bridge deck. The calculations yield an equivalent surface load of 2.491 kPa for the T structure area and 0.601 kPa at the span center, simulating the self-weight of the prestressed beams. Additionally, the permanent load of the bridge deck pavement is calculated at 4.801 kPa, while the variable load corresponds to the actual weight of the vehicles.

*3.4 Construction stages*

Using the 'Model Change' interaction, the concrete, rebars, and prestressed steel strands for each segment of the girder are activated in distinct analysis steps. Simultaneously, the corresponding self-weight loads and prestress tension for each segment are also activated. The activation sequence is as follows: (1) the 0# block (on the top of bridge pier) and side-span closure section; (2) the side span and mid-span segments (symmetric activation); and (3) the mid-span closure segment. When activating the side-span closure segment, boundary conditions are applied to limit the vertical displacement of the bottom plate in advance. These vertical constraints are then released once the side span is closed. According to the cooling method [23], the prestressed steel strands are stretched, with the coefficient of thermal expansion set to $2\times10^{-5}$ and the temperature difference defined as -370°C.

**4. Load test for FEM validation**

To verify the accuracy of the finite element model (FEM) and to investigate the stress mechanisms of the large-span, wide-box continuous rigid bridge, a load test was conducted on the actual structure. The measured deflection and strain at the control sections can be compared with the results obtained from the finite element analysis.

*4.1 Loading vehicles*

The load was applied by arranging four trucks for each section, with each truck weighing between 32.2 and 34.1 tons. To meet the specific requirements of the test, three-axle trucks were selected for loading. The loading test included both symmetrical and eccentric loading configurations along the cross-section. During the loading process, vehicles were loaded sequentially, starting from the smallest pile number to the largest, and from left to right across the cross-section. The loading sections included the mid-span (at 1/2 span), 1/4 of the mid-span, and both 1/3 and 1/2 of the side span. The loading was conducted in three stages, with an additional four vehicles added near the same section at each stage.

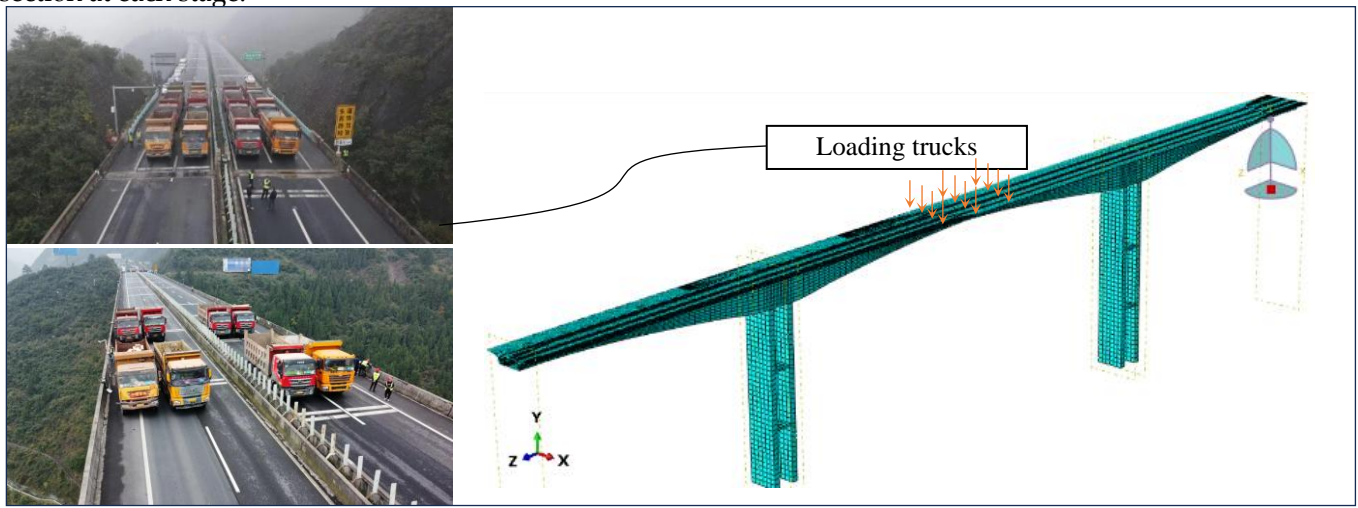

**Figure 3.** Load test and finite element model

*4.2 FEM validation*

The measured deflection values along the girder are compared with that calculated by the finite element modelling (FEM) in Figure 4. The FEM deflection under each loading stage aligns well with the measured site deflection in the loading test. It indicates that the finite element model is capable of accurately simulating the overall deformation of the bridge.

The strain results obtained from the finite element model (FEM) are compared with those from the loading test in Figure 5. The comparison shows that the FEM, which incorporates the construction process, accurately predicts the local strain distribution pattern of the bridge. However, it is important to note that the main girder exhibited significant downward deflection at mid-span and extensive cracking during service, caused by the excessive design aspect ratio and insufficient pre-stressing strand skew bending. Consequently, the strain gauges installed on-site may have been positioned over concealed cracks or between two cracks, leading to substantial discrepancies between the measured strain data and the FEM predictions. This issue is particularly evident at the mid-span, where the deflection is most pronounced and transverse cracking is most severe. The comparative results also reveal significant stiffness degradation near the middle of the span, further highlighting the critical structural deterioration in this region.

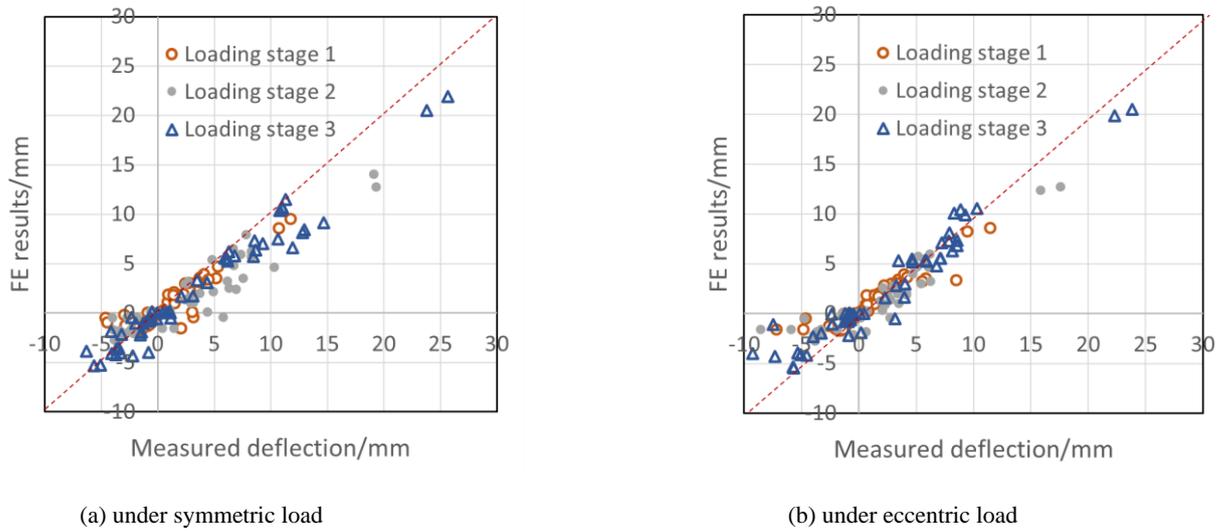

(a) under symmetric load

(b) under eccentric load

**Figure 4.** Comparison between FE analysis and test

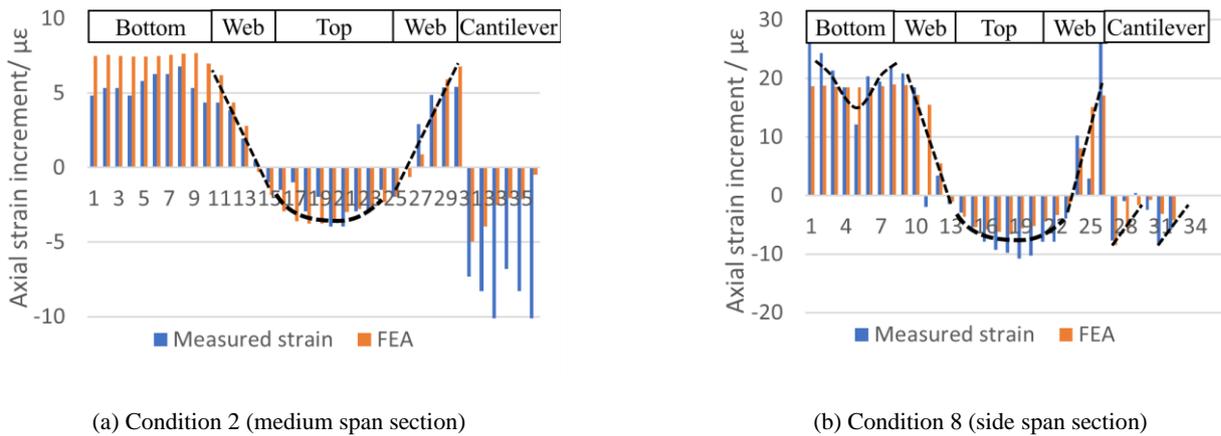

(a) Condition 2 (medium span section)

(b) Condition 8 (side span section)

Note: the abscissa is the section strain measuring point number, 1-10 is the bottom plate, 15-25 is the roof plate, 11-15,25-30 is the web plate, 31-36 is the cantilever plate

**Figure 5.** Comparison of cross-sectional strain increment

## 5. Shear-lag effect of single-box bridge girder

*5.1 Stress analysis*

Zhang Y. et al [26] derived the expression for the longitudinal stress at any spatial location of a simply supported box girder under concentrated load:

$$\sigma(x,y,z) = \frac{Pyz}{2I} + 2Ek_1Pl\omega(x,y)\left(\alpha_1 \sin\frac{\pi z}{l} - 9\alpha_2 \sin\frac{3\pi z}{l} + 25\alpha_3 \sin\frac{5\pi z}{l} - 49\alpha_4 \sin\frac{7\pi z}{l}\right) \quad (1)$$

A formulation for the longitudinal stress at any given spatial location subjected to a uniform load is presented as follows:

$$\sigma(x,y,z) = \frac{qyz}{2I}(l-z) - \frac{Ek_1q\omega(x,y)}{k^2}\left[\cosh(kz) - 1 + \frac{1-\cosh(kl)}{\sinh(kl)}\sinh(kz)\right] \quad (2)$$

However, precisely predicting the stresses through theoretical calculations alone is challenging due to the complex spatial distribution within the solid concrete. The boundary conditions of the main girder in the actual bridge are more complex than those assumed in simple theoretical models. At mid-span, the axial degrees of freedom are not fully constrained but are instead influenced by a spring boundary corresponding to the flexural stiffness at the abutment face. In the side span, the main girder is supported by a simple support at one end and a spring boundary at the other. These complex boundary conditions pose significant challenges for accurate calculations using conventional theoretical formulas. As a result, there has been an increasing reliance on the finite element method (FEM) for analyzing the spatial effects of the structure.

Currently, finite element analysis (FEA) is commonly used to obtain these predictions. Figure 6 presents the longitudinal stress simulation results for the bridge box girder, derived using a complex spatial solid finite element model. The results show that most of the stresses are concentrated at the flange near the web, with longitudinal stresses at the top plate near the web exceeding 15.0 MPa. In contrast, the stresses at the cantilever section, which is distant from the web, as well as at the flange near the centerline, are comparatively lower, with longitudinal stresses at the cantilever end and at the center of the top plate ranging from 7.3 to 8.2 MPa. The section exhibits a pronounced positive shear lag effect. Both the top and bottom plates show significant stress concentration near the center of the mid-span, and the uneven distribution of longitudinal stresses further highlights the shear lag effect. It is important to note that the diaphragm effectively restrains the mid-span section, leading to a more uniform distribution of stresses across both the top and bottom plates, and reducing stress concentration at the web-wing junction. However, sections farther from the diaphragm restraint experience a marked shear lag effect, which could result in insufficient pre-existing compressive stresses. As a result, if the variable load from an upper vehicle exceeds a certain threshold, it could cause concrete cracking, potentially compromising the durability and safety of the bridge.

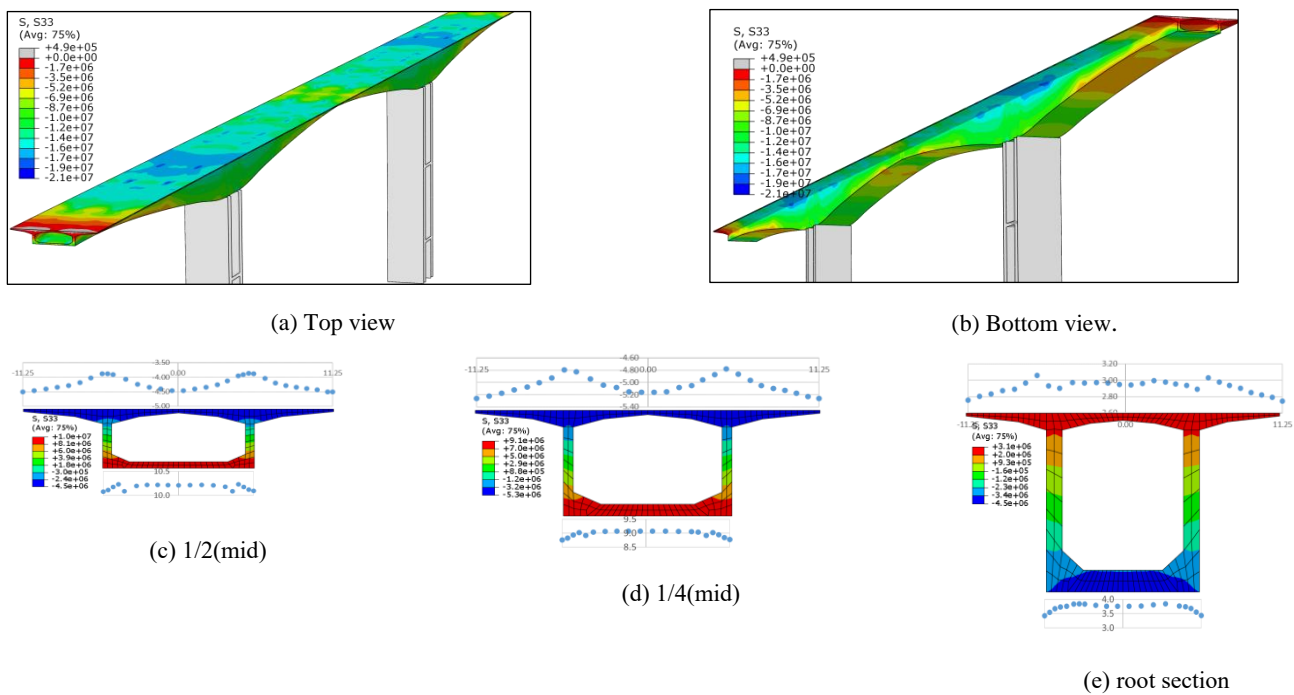

**Figure 6.** Longitudinal stress distribution in girder completed closure.

The main girder, especially for the large cantilever box girder, is subject to significant transverse bending, particularly in its upper flange, which may lead to longitudinal cracking on the road pavement, or even in the structural concrete. Figure 7 shows the transverse stress distribution in the main girder, highlighting that the transverse stresses at the top surface near the web are notably high, increasing the risk of transverse cracking. In contrast, the lower edge of the top plate's centerline experiences substantial transverse tensile stress, with stress levels approaching the ultimate tensile strength in the absence of transverse prestressing tendons. This condition makes the center of the internal top surface of the main girder especially susceptible to longitudinal cracking. Introducing transverse restraint at the bulkhead increases transverse tensile stress at its edge. Therefore, the application of transverse prestressing is essential to mitigate the risk of longitudinal cracking in the flange.

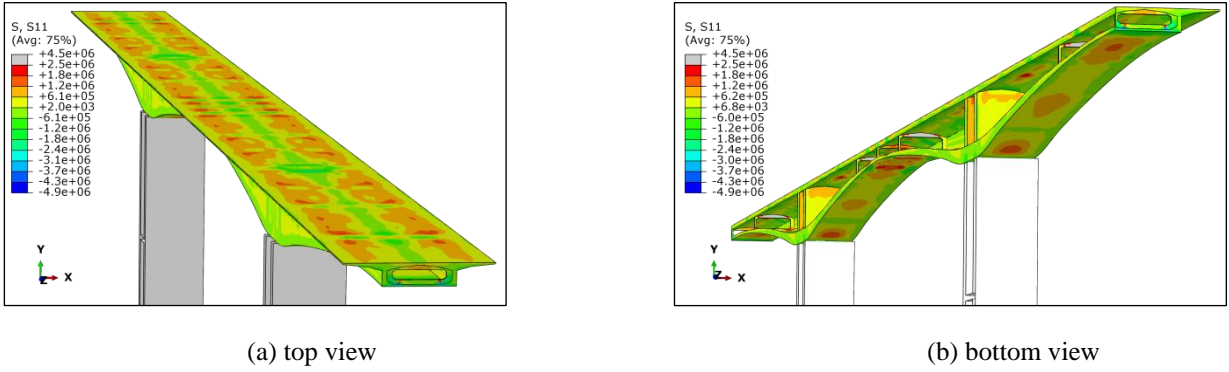

(a) top view  (b) bottom view

**Figure 7.** Transverse stress distribution of main girder after completion of the bridge

*5.2 Shear lag coefficient analysis*

The shear lag effect in the box structure arises from the non-uniform lateral stiffness along the transverse direction of the flange [27]. In regions adjacent to the web plate, both the shear stiffness and flexural stiffness increase, resulting in reduced bending deformation and increased positive bending stress. Conversely, in areas further from the web plate, the relative shear stiffness decreases, leading to greater bending deformation and a reduction in positive bending stress. Consequently, the longitudinal stresses are not uniformly distributed across the transverse direction [28]. As illustrated in Figure 8, this non-uniform distribution of longitudinal stress contributes to an uneven distribution of the shear lag coefficient along the transverse section. Notably, at the junction of the flange and the web, the shear lag coefficient deviates significantly from 1.0, indicating a pronounced shear lag effect in this region. The shear lag coefficient is calculated using Eq (3).

$$\lambda = \frac{\sigma}{\bar{\sigma}} \tag{3}$$

where $\lambda$, $\sigma$, $\bar{\sigma}$ are respectively the shear lag coefficient, the longitudinal normal stress of the bridge regard the spatial effect, and the normal stress of the longitudinal bridge calculated theoretically by the primary beam.

5.2.1. Shear lag coefficient under permanent load

To reduce self-weight, the main girder of a continuous rigid bridge typically employs a variable cross-section height, with both the flange and web plate thickness varying along the axial direction. Additionally, the plate thickness differs at various positions within the same cross-section. As a result, the spatial shear lag effect is highly complex and cannot be accurately captured using conventional theoretical formulas [29]. In this study, the width-to-height ratio of the main bridge box girder is 6.43, classifying it as a typical wide-box and large cantilever bridge. The shear lag effect, which results from stress non-uniformity, is particularly pronounced in this case. In the mid-span section, the spatial shear lag effect is particularly pronounced due to the sectional width-to-heigh ratio of 6.43, with the maximum shear lag coefficients of the top and bottom slabs being 1.09 and 1.56, respectively. It means the compressive stress in the concrete of this section may be up to 56% lower than the practical value calculated by primary beam theory.

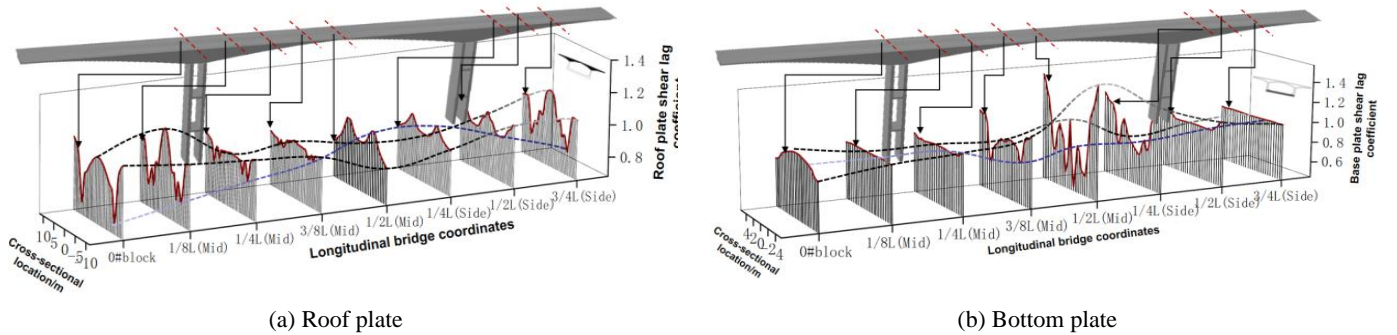

(a) Roof plate  (b) Bottom plate

**Figure 8.** Distribution of shear lag coefficient under permanent load.

To minimize the influence caused by the constraints of the bridge pier, a section located 5 meters away from the abutment edge was selected to be discussed. In regions with a large web height and negative bending moment, such as the 0# and 1/4 sections, the shear lag effect is typically expected to be minimal. However, the actual observation reveals that the shear lag coefficient at the center line of the top plate is significantly high, while the effect is less pronounced at other locations. The phenomenon is primarily attributed to the concentrated prestress at the junction of the web and flange. Consequently, this concentration of prestress counteracts the tensile stresses induced by the negative bending moment, resulting in lower stress and shear lag coefficients at the web-flange junction. A similar phenomenon can also be observed under live load conditions, as shown in Figure 9.

5.2.2. Shear lag coefficient under variable loads

To further examine the impact of vehicle load on shear lag, a finite element model of the bridge structure was developed, considering only concentrated loads. The results shown in Figure 9 demonstrate that the shear lag effect under variable concentrated loads follows a similar pattern to that observed under uniform loads. The maximum shear lag coefficient reaches 1.51, while the minimum is 0.64. As the applied load approaches the mid-span, the distribution of the shear lag coefficient becomes more uneven. Specifically, the maximum shear lag coefficient at each section increases gradually, peaking at the mid-span.

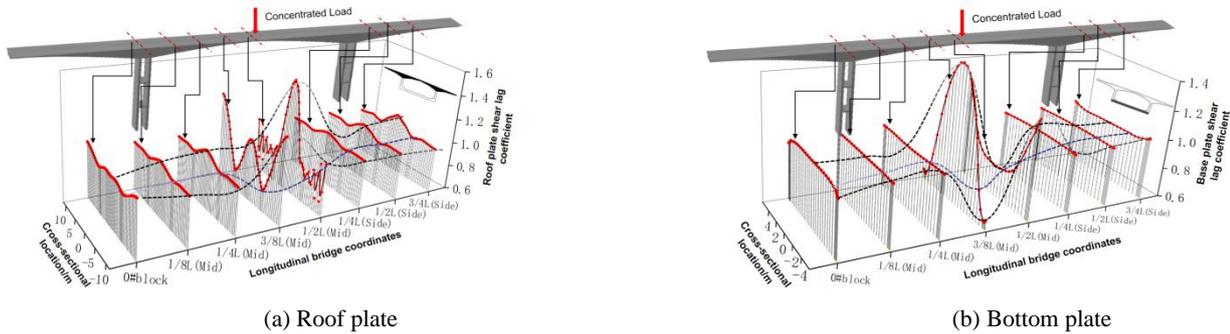

(a) Roof plate  (b) Bottom plate

**Figure 9.** Distribution of shear lag coefficients of the main beam under concentrated live load.

## 6. Parametric analysis

Previous studies have shown that the shear lag coefficient is significantly affected by the pair of box girder flange stiffness ratio $I_s/I$ [30], and the magnitude of the flange stiffness ratio of box girder is mainly determined by the section aspect ratio. When $I_s/I$ increases, i.e., the larger the section aspect ratio is, the larger the maximum shear lag coefficient of the section is. It can be seen that the larger the stiffness ratio of the flange plate is, the more serious the effect of shear lag is. For continuous rigid bridges with variable section height, the section aspect ratio at mid-span is the largest, so the shear lag coefficient is the largest and the shear lag effect is the most significant. In addition, the sensitivity factors such as width-to-span ratio, height-to-span ratio, width-to-thickness ratio, and wing plate width have equally important influences on the shear lag effect of variable-section box girders [31-36].

To this end, the models of continuous rigid bridge with variable cross-section main girder with different width-to-span and height-to-span ratios are established to analyze the variation rule of shear lag coefficient at different key cross-sections, including 1/8, 1/4 and 1/2 of mid span, and side-span 1/4-span and 1/2-span.

The main span of the bridge model is maintained at 210 m, with base plate half-widths of 5.25 m, 5.6 m, 5.96 m, and 6.3 m, corresponding to aspect ratios (2b/h) of 3, 3.2, 3.4, and 3.6, respectively. The conditions are simulated across different aspect ratios to observe the variations in the shear lag coefficients of the top and bottom plates. The results, illustrated in Figure 10, indicate that an increase in the overall section aspect ratio leads to a deviation of the shear lag coefficient from 1.0. This deviation arises because the change in section width enhances the full section stiffness $(I_s/I)$ of the upper and lower wing plates, thereby making the shear force lag effect more pronounced. Consequently, the maximum shear lag coefficient of the main girder section increases from 1.17 to 1.23, suggesting that the shear lag effect is more significant in wide box girder bridges compared to narrow box girder bridges.

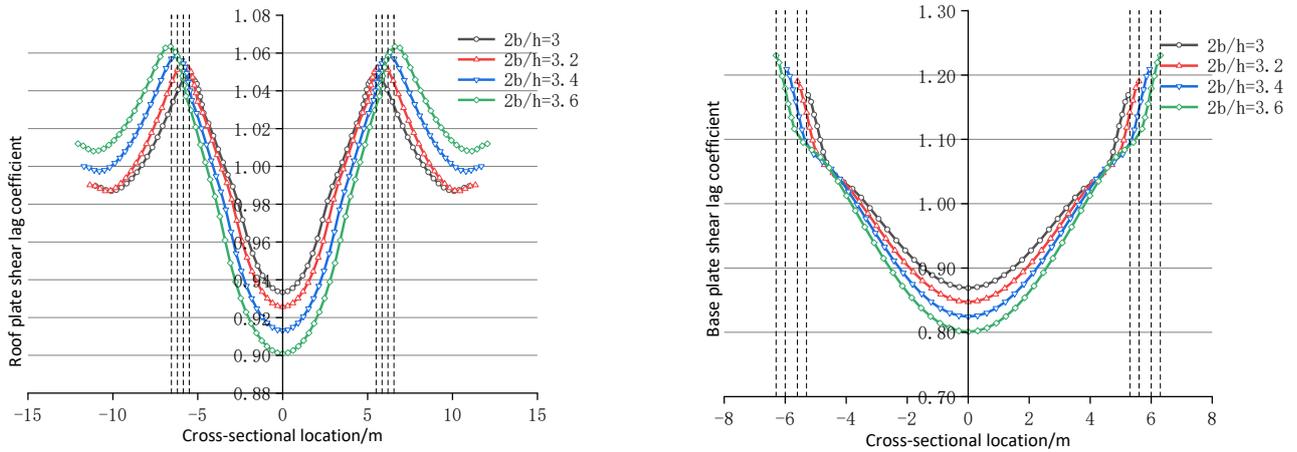

(a) 1/2(mid) span cross-section shear lag coefficient

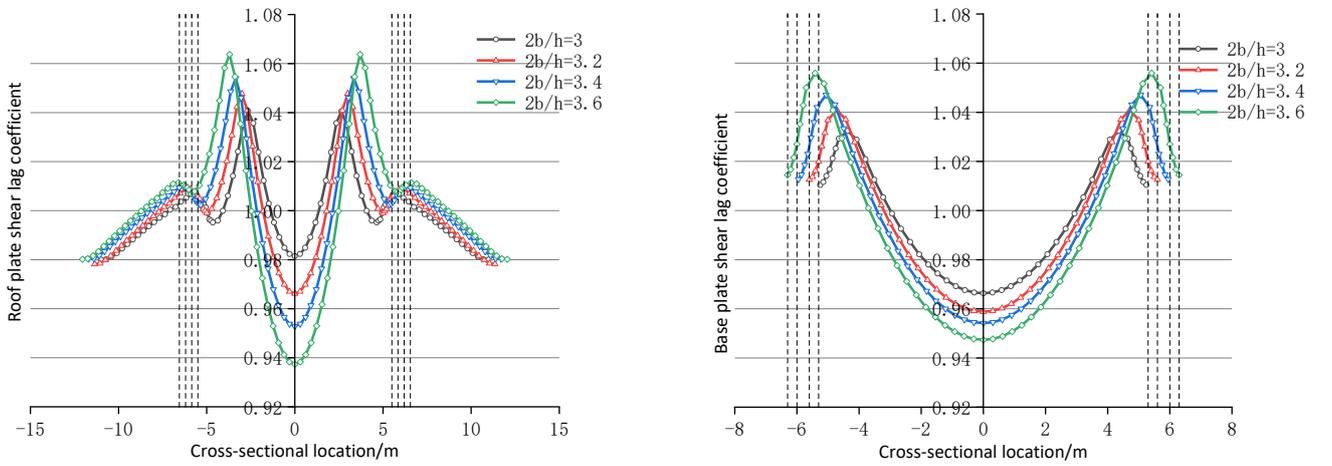

(b) 1/4(mid) span cross-section shear lag coefficient

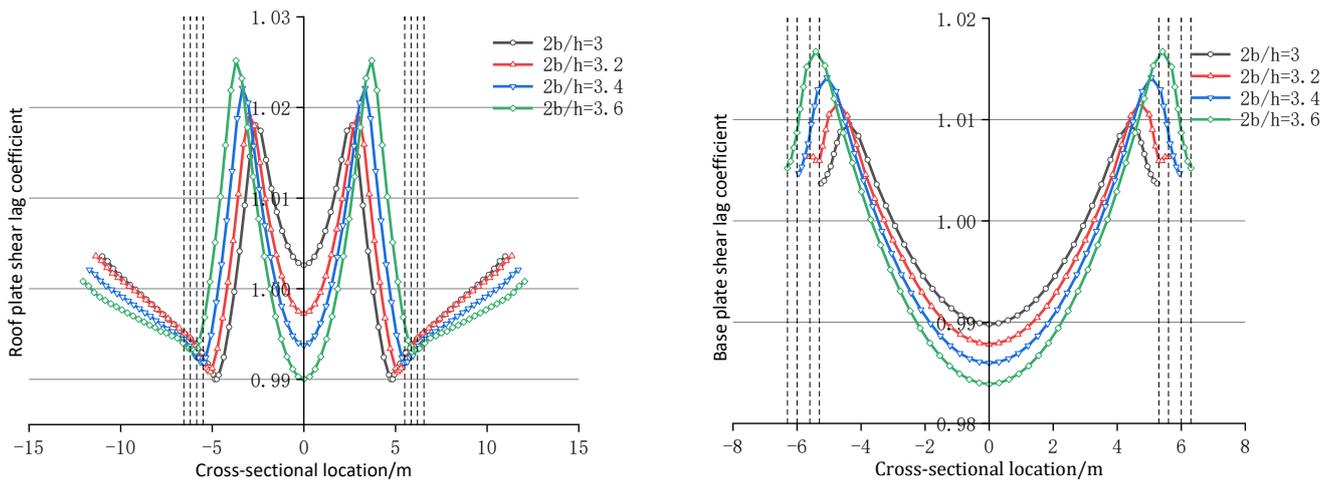

(c) 1/8(mid) span cross-section shear lag coefficient

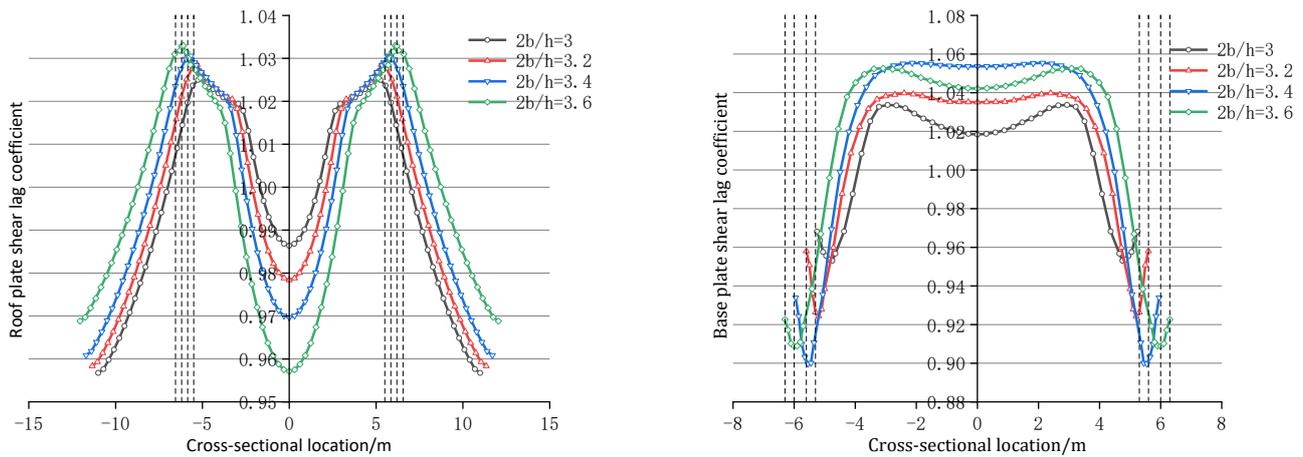

(d) 0#block cross-section shear lag coefficient

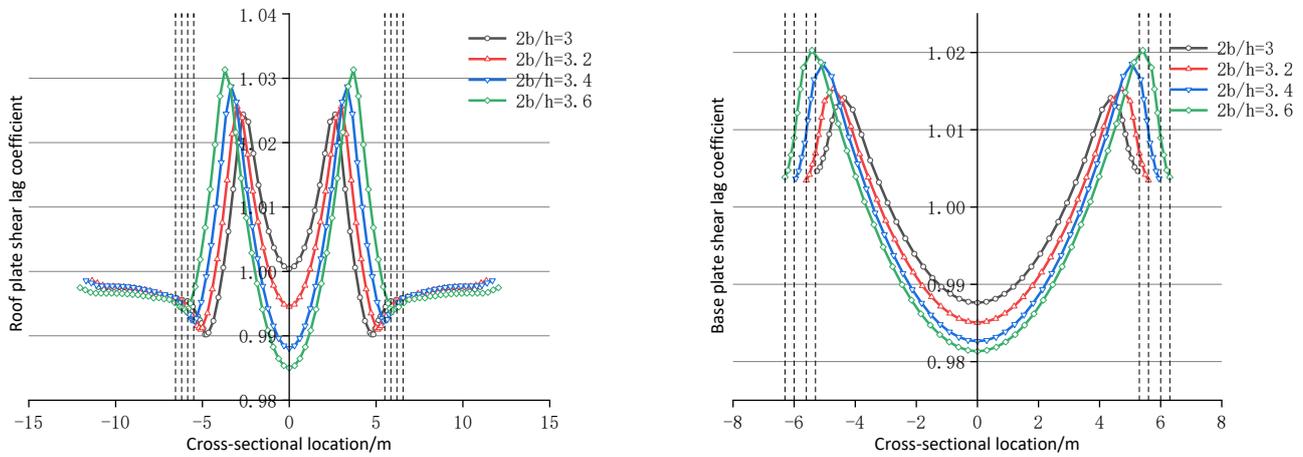

(e) 1/4(side) span cross-section shear lag coefficient

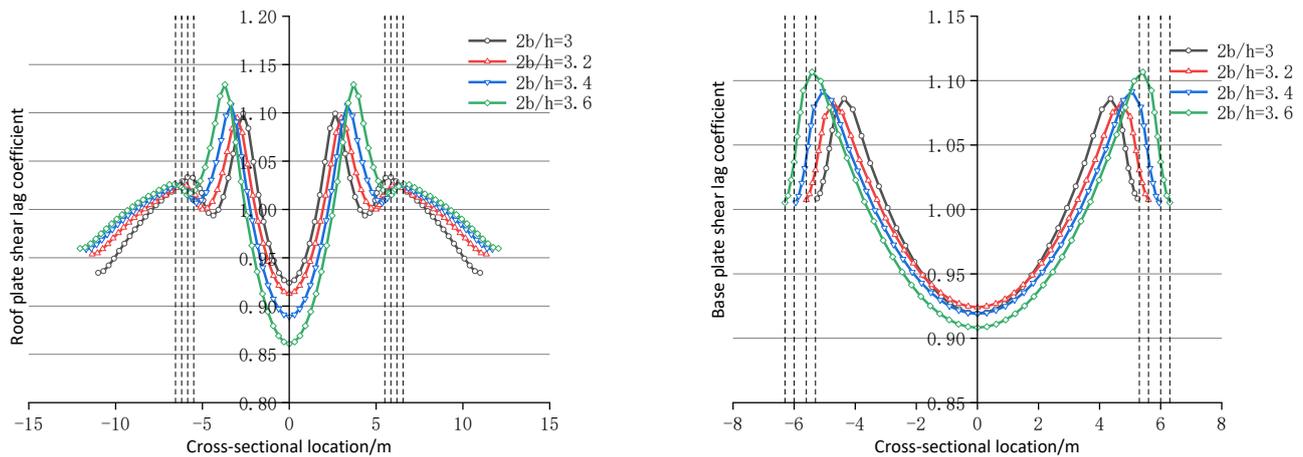

(f) 1/2(side) span cross-section shear lag coefficient

**Figure 10.** Distribution of shear lag coefficient at different aspect ratios

As the bridge span increases, the dead load effect exhibits a nonlinear increasing pattern. Under consistent width conditions, the main spans measure 165m, 185m, 210m, and 230m, with corresponding span-to-width ratios ($L/2b$) of 15, 17, 19, and 21, respectively. The variations in the shear lag coefficient of the roof and bottom plates at different span-to-width ratios are illustrated in Figure 11. Observing the bridge from a longitudinal perspective, it becomes evident that as the section moves further from the bearing position, the shear lag coefficient increasingly deviates from 1.0 at the same position along the section, with the mid-span shear lag effect being the most pronounced. This phenomenon is ultimately manifested in the maximum shear lag coefficient of the main beam section, which decreases from 1.35 to 1.16 as the span lengthens.

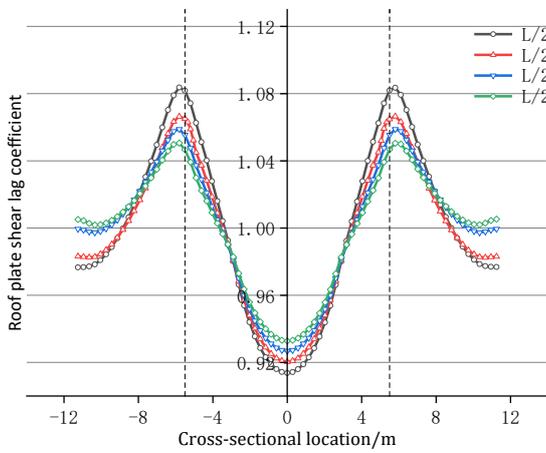 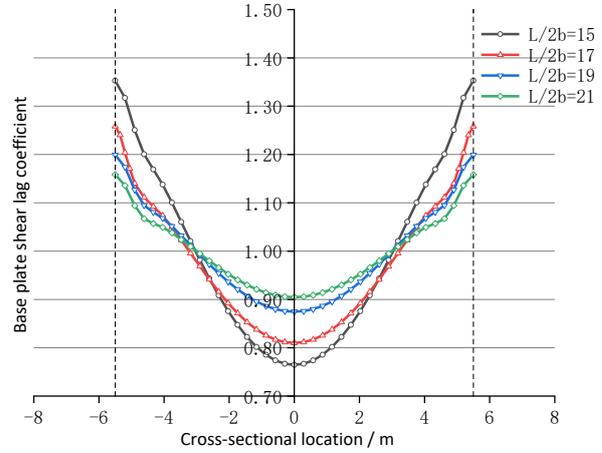

(a) 1/2(mid) span cross-section shear lag coefficient

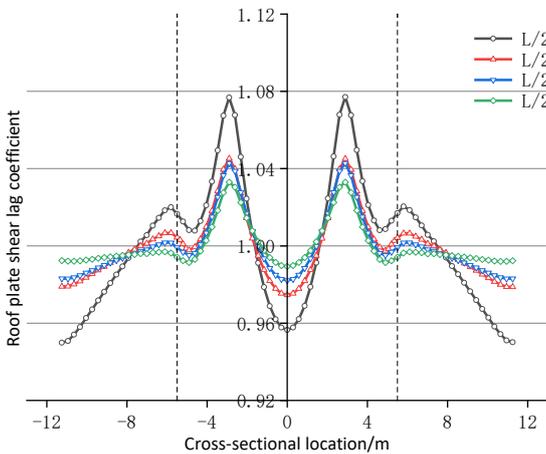 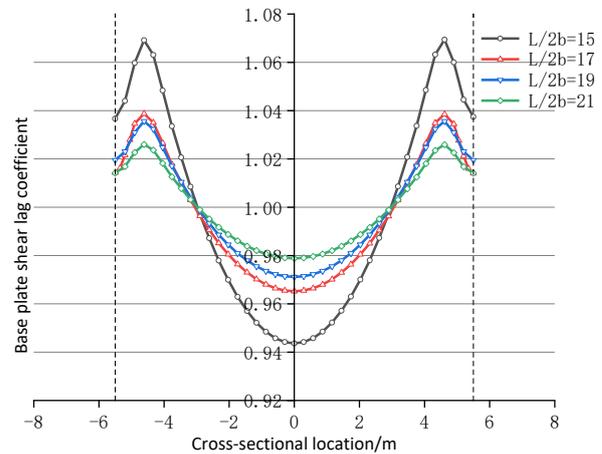

(b) 1/4(mid) span cross-section shear lag coefficient

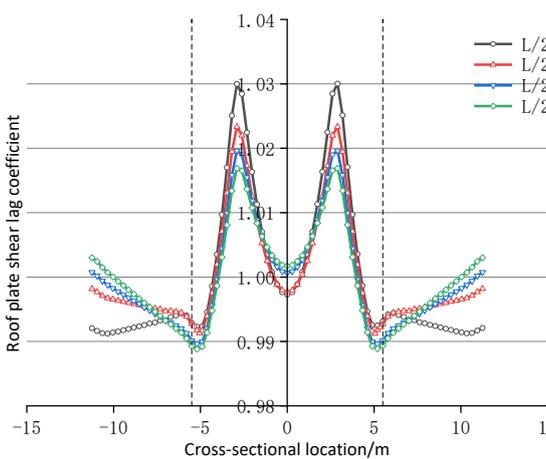 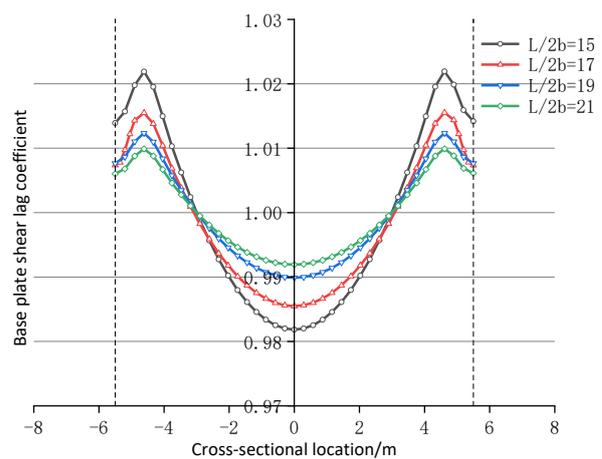

(c) 1/8(mid) span cross-section shear lag coefficient

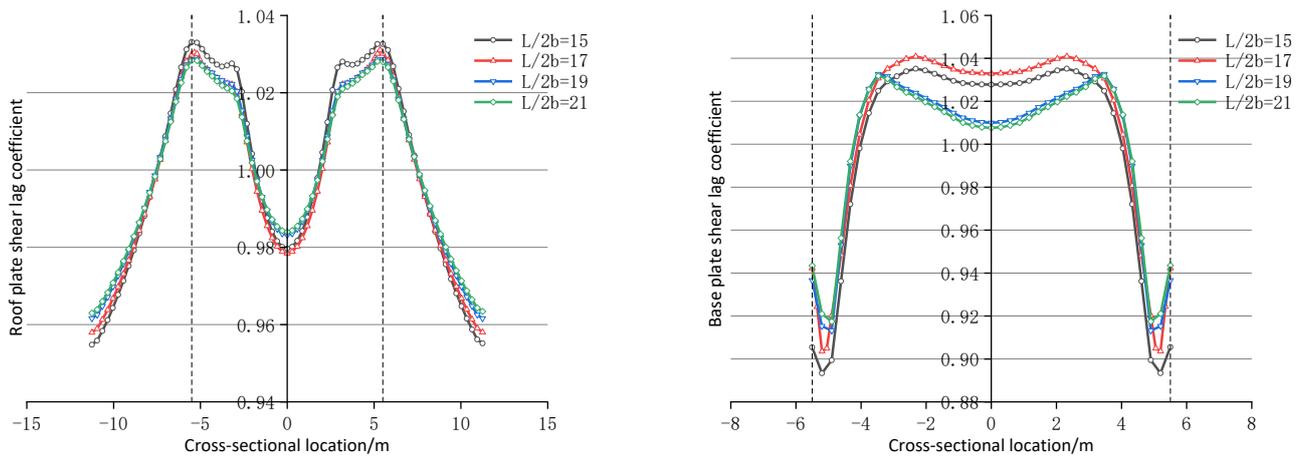

(d) 0# block cross-section shear lag coefficient

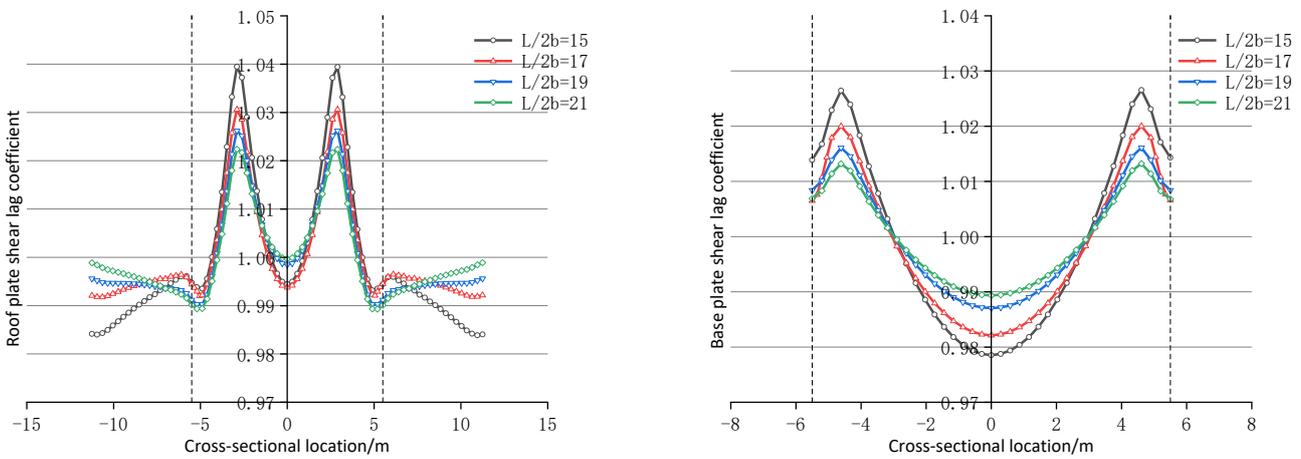

(e) 1/4(side) span cross-section shear lag coefficient

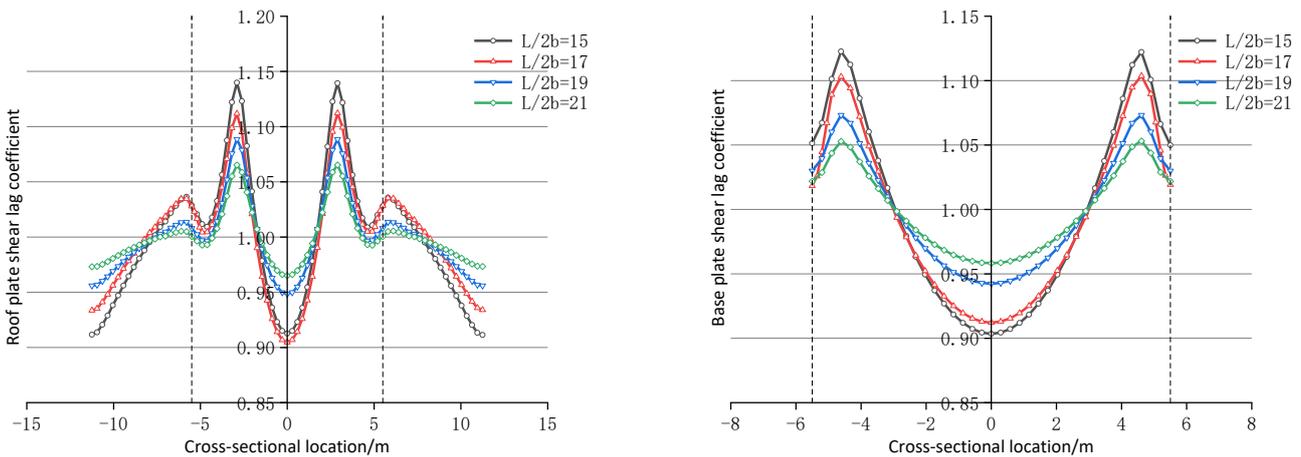

(f) 1/2(side) span cross-section shear lag coefficient

**Figure 11.** Distribution of shear lag coefficient at different span-to-width ratios

In general, the shear lag coefficient of the box girder exhibits a non-uniform distribution across the cross-section, affecting both the top and bottom plates. Due to the unique characteristics of the support and the influence of the relationship between the pier and the restraint, the shear lag effect at the top plate position tends to decrease with an increasing span ratio ($L/2b$) or an increasing aspect ratio ($2b/h$), and this effect is more remarkable. Yet, the bottom plate does not follow this trend, instead displaying an opposite behavior.

In summary, as the relative width of the flange increases, the span-to-width ratio (*L/2b*) decreases, while the aspect ratio (*2b/h*) increases. This results in a more pronounced spatial shear lag effect in the box girder, leading to an increase in longitudinal stress concentration at the intersection of the web and flange plates. Therefore, when designing long-span and wide box girder bridges, it is essential to fully consider the impact of the shear lag effect on the structural response to ensure the rationality and safety of the design.

**7. Conclusion**

To investigate the spatial shear lag effect of the large-span, wide-box prestressed concrete continuous rigid-frame bridge, a three-dimensional solid finite element model was developed, and load tests on a real bridge were performed to compare the deflection and local strain of the main girder to validate the model. In addition, the distribution of the shear lag coefficient in the main girder was analyzed under both uniform permanent and variable concentrated loads. After that, the variation of the shear lag coefficient with respect to different parameter conditions was also discussed. The main findings are as follows:

(1) Under different loading conditions, the shear lag coefficient at the flange-web junction exhibits a complex, non-uniform distribution both transversely and longitudinally, with significant deviations from 1.0. As the section approaches mid-span, the shear lag effect becomes more remarkable.
(2) At mid-span, where the width-to-height ratio is maximal, the shear lag effect is particularly significant. The maximum shear lag coefficient can reach to 1.51. In this region of positive bending moment, the pre-existing compressive stress in the bottom slab is substantially lower than that predicted by conventional beam theory, which may cause the stress at the bottom edge of the box girder to exceed the tensile strength of the concrete, potentially leading to cracking.
(3) In the negative bending moment region of the prestressed main girder, where the effects of self-weight and pre-stress partially counterbalance each other, the distribution of the shear lag coefficient deviates from the expected trend, with smaller shear lag coefficients observed at the flange-web junction.
(4) The parameter analysis reveals that the shear lag coefficients at the support section of the box girder differ from those at other sections, depending on the span-to-width ratio (*L/2b*) and the width-to-height ratio (*2b/h*). The variation is primarily attributed to the unique constraints and boundary conditions at the support locations.
(5) In the design of large-span wide box girders, the influence of shear lag effects should not be underestimated. Sensitivity analysis shows that as the span-to-width ratio (*L/2b*) decreases and the width-to-height ratio (*2b/h*) increases, the shear lag coefficient decreases from 1.35 to 1.16. Conversely, an increase in the width-to-height ratio results isn a rise in the shear lag coefficient from 1.17 to 1.23. These results demonstrate that both a smaller span-to-width ratio and a larger width-to-height ratio significantly affect the shear lag effect.

**Author Contributions:** Conceptualization, Zhaokun Shen; methodology, Ligui Yang and Shaorui Wang; software, Shuangshuang Jin and Chengmao He; validation, Chengtao Zhou; formal analysis, Dianhao Li; investigation, Lvye Zhou; resources, Yun Meng; data curation, Zhaokun Shen; writing—original draft preparation, Zhaokun Shen and Ligui Yang; writing—review and editing, Shaorui Wang ,and Shuangshuang Jin; visualization, Chengmao He; supervision, Yun Meng; project administration, Chengtao Zhou. All authors have read and agreed to the published version of the manuscript.

**Funding:** This research was funded by Project Supported by Department of Transport of Guizhou Province, grant number 2008009.

**Data Availability Statement:** All data included in this study are available upon request by contact with the corresponding author.

**Acknowledgments:** This research work is supported by the Project Supported by Department of Transport of Guizhou Province (2008009). The financial support is gratefully acknowledged

**Conflicts of Interest:** The authors declare no conflicts of interest.

**Reference**
1. Li Y.F., Xie J.Y., Wang F.C. & Li Y.H. Shear Lag Effect of Ultra-Wide Box Girder under Influence of Shear Deformation. *Special Issue-Bridge Structural Analysis*, *Applied Sciences*, 2024, 14(11): 4778.
2. Lian Y.D., Ma Q., Liu Z.Y., Zhang Y.Q. & Zhai Z.P. Shear Lag Effect of Framed Tube Structure under Horizontal Seismic Action. *Journal of Northwestern Polytechnical University*, 2017, 35(5): 890-897.
3. Liu F., Zhang Q. & Yang T. Effect of aspect ratio on shear distribution in box girders. *Engineering Structures*, 144: 12-22.
4. Huang M., LI Z. Research on Structural Performance of Prestressed Continuous Rigid Plate Bridge with Different Construction Sequences. *Highway*, 2024, 69 (02): 49-52.

5. Zhou J.S., Lou Z.H. The status quo and developing trends of large-span prestressed concrete bridges with continuous rigid frame structure. *China Journal of Highway and Transport*, 2000, 13(1): 31-37.
6. Zhou S.J., Jiang Y., Wu Y.D. Shear lag effect of twin-cell composite beams considering slip and shear deformation. *Civil Construction and Environmental Engineering*, 2017, 39(3): 20-27.
7. Zhou M.D., Li L.Y., Zhang Y.H. Research on Shear-lag Warping Displacement Function of Thin-walled Box Girders. *China Journal of Highway and Transport*, 2015, 28(6): 67-73.
8. Kármán T.V. Die mittragende breite. *Beiträge zur Technischen Mechanik und Technischen Physik*. 1924, 114-127.
9. Evans H.R., Ahmad M., Kristek V. Shear lag in composite box girders of complex cross-sections. *Journal of Constructional Steel Research*, 1993,24(3): 183-204.
10. Wang C.F., Zhu D.S., Chen X.C. Analysis of the shear lag effect of a double-chamber box girder bridge: The 12th National Conference on Structural Engineering, Chongqing, China, 2003.
11. Zhao Q.Y., Zhang Y.H., Shao J.Y. Analysis of shear lag effects in box girders with variable-thickness flanges.*Applied Mathematics and Mechanics*,2019,40(6): 609-619.
12. Tintner G. A Note on Welfare Economics. *Econometrica*, 1946,14(1): 69.
13. Li X.Y., Fan W., Wan S. Deflection calculation analyses on thin-walled box girder based on the theory of Timoshenko beam and the energy-variation principle. *Journal of South China University of Technology (Natural Science Edition)*, 2018,46(4): 8.
14. Hildebrand Francis B. The exact solution of shear-lag problems in flat panels and box beams assumed rigid in the transverse direction. *NASA Technical Notes,* 1943, No. 894.
15. De F.A, Scordelis A.C. Direct Stiffness Solution for Folded Plates. *Journal of the Structural Division*, 1964,90(4): 15-47.
16. Cai H., Lu H.L., Tang Z. Shear Lag Effect of Thin-walled Box Girder Based on Shell Finite Element Theory. *Journal of Railway Science and Engineering*, 2017,14(04): 779-786.
17. Wu Y.M., Luo Q., Yue Z.F. Finite element method for calculating the shear lag in box girders. *China Railway Science*, 2003, 24(4): 65∼69.
18. Zhang L., Xu H. & Wang H. Finite element analysis of shear hysteresis in continuous rigid bridges. *Computational Materials Science*, 2019, 155: 392-402.
19. Xu, J., Huang S. Shear flow and hysteresis effects in large-span box girders. *Journal of Bridge Engineering*, 25(4), 04020018.
20. Lin P.Z., Sun L.X., Yang Z.J. Research on Shear Lag Effect of Twin-cell Box Girders. *Journal of Railway Engineering Society*, 2014, 31(1): 59-63, 112.
21. Zhao L., Jiang F.X., Ma Z.Z. Effect of diaphragm setting on shear lag of concrete box girder under eccentric load. *East China Highway*, 2009(5): 46-48.
22. Li Y., Zhao Z. Numerical simulation of shear hysteresis in bridge box girders under dynamic loads. *Structural Engineering International*, 31(3), 346-358.
23. Yapar O., Basu P.K., N. Nordendale. Accurate finite element modeling of pretensioned prestressed concrete beams. *Engineering Structures*, 2015. 101:163-178.
24. GB 50010-2010. Code for Design of Concrete Structures. Beijing: China Architecture & Building Press, 2015.
25. Yang L.G., Wang Y.Y. Elchalakani M. and Fang Y. Experimental Behavior of Concrete-Filled Corrugated Steel Tubular Short Columns under Eccentric Compression and Non-uniform Confinement. *Engineering Structures* 2020, 220: 111009.
26. Zhang Y.Y., Zhang Y.H., Zhang H. A Separate Solution Method for Shear Lag Effects in Box Girders and Parameter Analysis. *Applied Mathematics and Mechanics*, 2018, 39(11):1282-1291.
27. Gao Y., Li Z.L. Shear Lag Effect in Tall Tubular Structures. *Journal of Southwest University of Science and Technology*, 2006, 21(2): 15-19.
28. Li H.H., Li L.F., Yoo D.Y., Ye M., Zhou C., Shao X.D. Experimental and theoretical investigation of the shear lag effect in the novel non-prismatic prestressed CSW–UHPC composite box girders. *Archives of Civil and Mechanical Engineering*, 2023, 23:132-154.
29. Chen H.M., Qiao J.Y., State of the art about shear lag effect of concrete box girder bridges. *Structural Engineers*, 2011, 27(1): 161-166.
30. Xiang H.F., Fan L.C. Advanced Theory of Bridge Structures. Beijing: China Communications Press, 2013: 14-19.
31. Cai H., Lu H.L., Tang Z. Vibration properties research on curved box girder considering shear lag effects. *World Earthquake Engineering*, 2016, 32(4): 239-244.
32. Liu Q.G. Shear Lag Effect of Variable Box Section Continuous Girder and Its Formation Mechanism. *Railway Construction Technology*, 2020, (04): 6-9.
33. Gao L., Zhou Y., Wang S.L. Energy-variational method of shear-lag effect in PC box-girder with corrugated steel webs. *Journal of Xi' an University of Architecture and Technology ( Natural Science Edition)*, 2022, 54 (1): 27-34.
34. Wu W.Q., Wan S., Ye J.S. 3-D finite element analysis on shear lag effect in composite box girder with corrugated steel web.*China Civil Engineering Journal*, 2004, 37(9): 31-36.
35. Shu F. M., Xie D.Z. Parametric Analysis on Shear Lag Effects of Box-Girder Beam. *Applied Mechanics and Materials*, 2013, 2544(351352): 152-155.
36. Liu X.Z., Cheng K., Chen Y. Shear lag effect of long span composite continuous box-girder bridge with corrugated steel webs and variable cross-section. *Journal of Shenyang jianzhu university ( nature science)* , 2017, 33(5): 855-862.

37. Zhang Y.H. The experimental study and theoretical analysis in shear lag effect on thin-walled box-girder. Lanzhou: Lanzhou Jiaotong University, 2011.
38. Zhou S.J. Shear lag analysis of box girders. *Engineering Mechanics*,2008, 25(2):204-208.
39. Ding W.L. Analysis of shear Lag Effect of Large-span Prestressed Concrete Continuous Box Girder Bridge .Changsha: Central South University,2009.
40. Sun Z.W., Liu Y., Wu X.G., et al. Analysis of shear lag effect and effective distribution width of continuous rigid frame wide box girder. *Journal of Railway Science and Engineering*, 2016, 13(7): 1347–1351.